\documentclass{ifacconf}
\usepackage{graphicx}      
\usepackage{natbib}        
\def\disp{\displaystyle}
\def\dref#1{(\ref{#1})}
\def\crr{\cr\noalign{\vskip1.5mm}}
\usepackage{mathrsfs}
\usepackage{amsfonts}
\usepackage{amsmath,amssymb}
\newtheorem{theorem}{Theorem}[section]   

\newtheorem{lemma}{Lemma}[section]

\newcommand{\rline}  {{\mathbb R}}

\renewcommand{\e}    {{\varepsilon}}

\newcommand{\m}      {{\hbox{\hskip 1pt}}}

\newcommand{\dd}     {{\rm d\hbox{\hskip 0.5pt}}}
\newcommand{\Ascr}   {{\cal{A}}}
\newcommand{\Bscr}   {{\cal{B}}}

\begin{document}
\begin{frontmatter}

\title{Output feedback exponential stabilization of a nonlinear 1-D
       wave equation with boundary input \thanksref{footnoteinfo}}

\thanks[footnoteinfo]{This work was partially supported by grant no.
        800/14 of the Israel Science Foundation.}

\author[First]{Hua-Cheng Zhou}
\author[First]{George Weiss}

\address[First]{School of Electrical Engineering, Tel Aviv University,
 Ramat Aviv, Israel (e-mail: hczhou@amss.ac.cn)
 (gweiss@eng.tau.ac.il).}

\begin{abstract}
This paper develops systematically the output feedback exponential
stabilization for a one-dimensional unstable/anti-stable wave equation
where the control boundary suffers from both internal nonlinear
uncertainty and external disturbance. Using only two displacement
signals, we propose a disturbance estimator that not only can estimate
successfully the disturbance in the sense that the error is in
$L^2(0,\infty)$ but also is free high-gain. With the estimated
disturbance, we design a state observer that is exponentially
convergent to the state of original system. An observer-based output
feedback stabilizing control law is proposed. The disturbance is then
canceled in the feedback loop by its approximated value. The
closed-loop system is shown to be exponentially stable and it can be
guaranteed that all internal signals are uniformly bounded.
\end{abstract}

\end{frontmatter}

\section{Introduction} 

We consider the output feedback exponential stabilization problem of
a one-dimensional unstable nonlinear wave equation with boundary
input $u+d$:
\begin{equation} \label{wave-o}
   \left\{\begin{array}{l}
   w_{tt}(x,t) \m=\m w_{xx}(x,t),\ \ 0<x<1,\ t>0, \crr\disp
   w_x(0,t) \m=\m -qw(0,t),\ \ t\geq0, \crr\disp
   w_x(1,t) \m=\m u(t)\hspace{-0.1cm}+\hspace{-0.1cm}f(w(\cdot,t),
   w_t(\cdot,t))\hspace{-0.05cm}+\hspace{-0.05cm}d(t),\ t\geq0,
      \crr\disp
   w(x,0) \m=\m w_0(x),\;w_t(x,0)=w_1(x),\;0\leq x\leq1,\crr\disp
   y_m(t) \m=\m \{w(0,t),w(1,t)\}. \end{array}\right.
\end{equation}
Here $(w,w_t)$ is the state, $u$ is the control input signal, $y_m$ is
the output signal, that is, the boundary traces $w(0,t)$ and $w(1,t)$
are measured. The equation containing the constant $q>0$ creates a
destabilizing boundary feedback at $x=0$ that acts like spring with
negative spring constant. The function $f:H^1(0,1)\times
L^2(0,1)\to\rline$ is unknown and represents the internal uncertainty,
while $d$ represents the unknown external disturbance which is only
supposed to satisfy $d\in L^\infty(0,\infty)$. We use the notation
\begin{equation} \label{Helen}
   F(t) \m=\m f(w(\cdot,t),w_t (\cdot,t))+d(t)
\end{equation}
for the ``total disturbance''. We consider system \dref{wave-o} in the
state Hilbert space $\mathbb{H}=H^1(0,1)\times L^2(0,1)$ with the
usual inner product. Our aim is to design a feedback controller which
generates the control signal $u$ (using the measurements $y_m$) such
that the state of the system converges to zero, exponentially.

Later in the paper, we shall also discuss a related problem, where
the negative spring is replaced by a negative damper. More precisely,
on the right hand-side of the equation containing $q$, we have
$-qw_t(0,t)$. We shall solve the exponential stabilization problem
also for this alternative nonlinear wave system.

It is well known that output feedback stabilization is one of the
fundamental issues in control theory. The key idea in output feedback
is that the control and output should be as little as possible. When
the internal uncertainty and the external disturbance flow in the
control end, the stabilization problem \dref{wave-o} becomes much more
complicated.  In this paper, we present a dynamic compensator which
employs a PDE disturbance estimator and full state feedback based on
the observer state. Our compensator consists of two parts: the first
part is to cancel the total disturbance by applying the active
disturbance rejection control (ADRC) strategy, which is an
unconventional design strategy first proposed by Han in
\cite{Han1998}; the second part is to stabilize system by using the
backstepping approach. The stabilization problem of system
\dref{wave-o} was first considered in \cite{GWGBZ2013tac}, where the
output measurement is $y_m(t)=\{w(0,t),w_t(1,t)\}$, the adaptive
controller is designed, and the disturbance $d$ has the following form:
$d(t)=\sum_{j=1}^m[\bar{\theta}_j\sin\alpha_jt+\bar{\vartheta}_j\cos
\alpha_jt],\; t\geq0 $ with known frequencies $\alpha_j,$ and unknown
amplitudes $\bar{\theta}_j,\bar{\vartheta}_j,j=1,2,\ldots,m$, and the
resulting closed-loop system is asymptotically stable.  Obviously, the
disturbance signal in our paper is more general than the harmonic
disturbance. Recently, the stabilization problem of system
\dref{wave-o} with $f\equiv0$, $d\in L^\infty(0,\infty)$ has been
investigated in \cite{FG2016ccc}, where the output measurements are
$\{w(0,t),w_t(0,t),w(1,t)\}$, and their result is that the
closed-loop system is asymptotically stable.
The output feedback in \cite{FG2016ccc} uses one more measurement than
\cite{GWGBZ2013tac}.  
Another point that in our treatment is different from \cite{GWGBZ2013tac,FG2016ccc} is that the
closed-loop systems in our paper are exponential stable and we do not
require to measure the velocity $w_t(0,t)$ (or $w_t(1,t)$) which is
hard to measure, as explained in (\cite{FansonVmeas}). In this paper,
we only use two scalar signals (the components of $y_m$). It can be
shown that the problem can hardly be solved with less than two output
signals, which will be shown in the journal version of this work.

Output feedback stabilization for one-dimensional anti-stable wave
equation was considered in \cite{JFF2015tac}, where a new type of
observer is constructed by using three output signals to estimate the
state first and then estimate the disturbance via the state of
observer through an extended state observer (ESO). However, the
initial value is required to be smooth in \cite{JFF2015tac} and they
obtain asymptotic stability (not exponential, like here). In the
recent work \cite{FHGBZ}, the author introduces a new disturbance
estimator which is different from the traditional estimator, the
smoothness requirement on the initial state being removed. In
\cite{FHGBZ}, still three output signals are used and the controller
achieves asymptotic stability.

The paper is organized as follows: in Section 2 we design an
infinite-dimensional disturbance estimator that does not use
high gain. We propose a state observer based on this estimator and
we develop an output feedback controller in Section 3. The
exponential stability of the closed-loop system is proved in Section
4. Finally, Section 5 is devoted to the output feedback exponential
stabilization of the alternative anti-stable wave equation mentioned
earlier (with the negative damper).

\section{Disturbance estimator design}\label{Sec-dis-design} 

The following lemma is not difficult to prove by using the results in
\cite{Weiss.G} and \cite{TWeiss2009book}. For related results we refer
to \cite{Jacob2016}.

\begin{lemma} \label{Lem-ABu}
Let $A$ be the generator of exponential stable $C_0$-semigroup
$e^{At}$ on the Hilbert space $X$. Assume that $B_i\in
\mathcal{L}(U_i,X_{-1})$, $i=1,2,\ldots,n$ are admissible control
operators for $e^{At}$. Then, the initial value problem
$\dot{x}(t)=Ax(t)+\sum_{i=1}^nB_iu_i(t),\; x(0)=x_0,\;u_i\in
L^2_{loc}(0,\infty;U_i)$, admits a unique solution $x\in
C(0,\infty;X)$, which tends to zero as $t\to\infty$ if either $u_i\in
L^2(0,\infty;U_i)$ or $\lim_{t\to\infty}\|u(t)\|_{U_i}=0$, for
$i=1,2,\ldots,n$, and is bounded if $u_i\in L^\infty(0,\infty;U_i)$,
$i=1,2,\ldots,n$. Moreover, if there exist two constants $M_0,\mu_0>0$
such that $\|u\|_{U_i}\leq M_0e^{-\mu_0 t}$, $i=1,2,\ldots,n$, then
$\|x(t)\|\leq Me^{-\mu t}$ for some $M,\mu>0$.
\end{lemma}

Now we design a {\em total disturbance estimator} for the system
\dref{wave-o}. This is an infinite dimensional system with the state
consisting of the functions $v,v_t,z,z_t,W$, defined on $(0,1)$:
\begin{equation} \label{wave-o-v_first}
\left\{\begin{array}{l}
   v_{tt}(x,t) \m=\m v_{xx}(x,t),\ \ 0<x<1,\;t>0, \crr
   v_x(0,t) \m=\m -qw(0,t)+c_1[v(0,t)-w(0,t)], \crr
   v_x(1,t) \m=\m u(t)-W_x(1,t),\ \ t\geq0, \crr
   v(x,0) \m=\m v_0(x),\ \ v_t(x,0)=v_1(x),\ \ 0\leq x\leq 1,\crr
   W_t(x,t) \m=\m -W_x(x,t),\;0<x<1,\ \ t>0, \crr
   W(0,t) \m=\m -c_0[v(0,t)-w(0,t)],\ \ t>0, \crr
   W(x,0) \m=\m W_0(x),\ \ 0\leq x\leq 1,\end{array}\right.
\end{equation}
\begin{equation} \label{wave-o-v_second}
\left\{\begin{array}{l}
   z_{tt}(x,t) \m=\m z_{xx}(x,t),\;0<x<1,\;t>0, \crr
   z_x(0,t) \m=\m \frac{c_1}{1-c_0}z(0,t)+\frac{c_0}{1-c_0}
            z_t(0,t),\ \ t\geq 0, \crr
   z(1,t)=v(1,t)+W(1,t)-w(1,t),\ \ t\geq 0, \crr
   z(x,0)=z_0(x),\;z_t(x,0)=z_1(x),\ \ 0\leq x\leq 1.\end{array}\right.
\end{equation}
Here $c_0$ and $c_1$ are two positive design parameters,
$(v_0,v_1.z_0,z_1,W_0)\in \mathbb{H}^2\times H^1(0,1)$ is the initial
state of the disturbance estimator, and its inputs are $u$ and $y_m$.
The part \dref{wave-o-v_first} is used to channel the total disturbance
$F$ to an exponentially stable system. Indeed, set $\widehat{v}(x,t)=
v(x,t)-w(x,t)$. Then $(\widehat{v}(x,t),W(x,t))$ satisfies
\begin{equation} \label{wave-o-hatv-P}
   \left\{\begin{array}{l}
   \widehat{v}_{tt}(x,t)=\widehat{v}_{xx}(x,t),  \crr\disp
   \widehat{v}_x(0,t)=c_1\widehat{v}(0,t),  \;
   \widehat{v}_x(1,t)+W_x(1,t)=-F(t), \crr\disp
   W_t(x,t)=-W_x(x,t),  \;\; 
   W(0,t)=-c_0\widehat{v}(0,t). \end{array}\right.
\end{equation}
It follows from the next lemma that the linear part of
\dref{wave-o-hatv-P} (when $F=0$) is exponentially stable. We remark
that the well-posedness assumption about system \dref{wave-o} (which
appears in the lemma) can actually be proved, and this will be in the
journal version of this work.

\begin{lemma} \label{Lem-hatz-bd}
Suppose that $c_0\in(0,1)$, $d\in L^\infty(0,\infty)$ (or $d\in
L^2(0,\infty)$), $f:\mathbb{H}\to\mathbb{R}$ is continuous and that
\dref{wave-o} admits a unique solution $(w,\dot{w})\in
C(0,\infty;\mathbb{H})$ which is bounded. For any initial value
$(\widehat{v}_0,\widehat{v}_1,W_0)\in \mathbb{H}\times H^1(0,1)$ with
the compatibility condition $W_0(0)=-c_0\widehat{v}_0(0)$, the system
\dref{wave-o-hatv-P} admits a unique solution $(\widehat{v},
\widehat{v}_t,W)\in C(0,\infty;\mathbb{H}\times H^1(0,1))$ such that
$$ \sup_{t\geq0}\|(\widehat{v}(\cdot,t),\widehat{v}_t(\cdot,t),
   W(\cdot,t))\|_{\m\mathbb{H}\times H^1(0,1)} \m<\m+\infty.$$
Moreover, if $\lim_{t\to\infty}|f(w,w_t)|=0$ and $d\in L^2(0,\infty)$,
then $\lim_{t\to\infty}\|(\widehat{v}(\cdot,t),\widehat{v}_t(\cdot,t),
W(\cdot,t))\|_{\mathbb{H}\times H^1(0,1)}=0$. If $f\equiv 0$
and $d\equiv0$, then $\|(\widehat{v}(\cdot,t),\widehat{v}_t(\cdot,t),
W(\cdot,t))\|_{\mathbb{H}\times H^1(0,1)}\leq M'e^{-\mu' t}$ for
all $t\geq0$ with some $M',\mu'>0$.
\end{lemma}

{\bf Proof.} Let $\widetilde{v}(x,t)=\widehat{v}(x,t)+W(x,t)$, then
it is easy to check that $(\widetilde{v}(x,t),W(x,t))$ is governed by
\begin{equation}\label{wave-o-hatv}
\left\{\begin{array}{l}
\widetilde{v}_{tt}(x,t)=\widetilde{v}_{xx}(x,t),  \cr 
\widetilde{v}_x(0,t)=\frac{c_1}{1-c_0}\widetilde{v}(0,t)+\frac{c_0}
   {1-c_0}\widetilde{v}_t(0,t),  \cr \disp
\widetilde{v}_x(1,t)=-f(w(\cdot,t),w_t(\cdot,t))-d(t),  \cr \disp
W_t(x,t)=-W_x(x,t), \;\; 
W(0,t)=-\frac{c_0}{1-c_0}\widetilde{v}(0,t).
\end{array}\right.
\end{equation}
We first consider the ``$\widetilde{v}$-part'' of \dref{wave-o-hatv}.
To this end, define the operators $\Ascr$ and $\Bscr$ by:
$ \Ascr(\phi,\psi)^\top \m=\m (\psi,\phi'')^\top$, $\forall\ (\phi,
   \psi)^\top\in D(\Ascr)$ with
   $ D(\Ascr) \m=\m \big\{ (\phi,\psi)^\top\in H^2(0,1)\times
                  H^1(0,1)\ |\phi'(1)$ $=0,\phi'(0)=\frac{c_1}{1-c_0}\phi(0) +
   \frac{c_0}{1-c_0}\psi(0)\big\},$
and $\Bscr=(0,-\delta(x-1))^\top$. Then the ``$\widetilde{v}$-part'' of
\dref{wave-o-hatv} can be written  as
$$ \dfrac{\dd}{\dd t}\begin{pmatrix}\widetilde{v}(\cdot,t)\cr
   \widetilde{v}_t(\cdot,t)\end{pmatrix} = \Ascr\begin{pmatrix}
   \widetilde{v}(\cdot,t)\cr \widetilde{v}_t(\cdot,t)\end{pmatrix}
   + \Bscr F(t) \m,$$
where $F$ is given by \dref{Helen}. It is well-known that $\Ascr$
generates an exponential stable $C_0$-semigroup $e^{{\mathcal A}t}$
and $\mathcal{B}$ is admissible to $e^{{\mathcal A}t}$. Since
$f:\mathcal{B}\to\mathbb{R}$ is continuous and $(w,\dot{w})^\top\in
C(0,\infty;\mathbb{H})$ is bounded, $f(w)\in
L^\infty(0,\infty)$. Thus, by $d\in L^\infty(0,\infty)$ or by $d\in
L^2(0,\infty)$, it follows from Lemma \ref{Lem-ABu} that
``$\widetilde{v}$-part'' of \dref{wave-o-hatv} admits a unique bounded
solution, i.e., $\sup_{t\geq0}\|(\widetilde{v}(\cdot,t),\widetilde
{v}_t(\cdot,t))^\top\|_{\mathbb{H}}\leq M_1 $ with $M_1>0$.
Next, we claim that $\sup_{t\geq0}\|W(\cdot,t))\|_{H^1(0,1)}<+\infty$.
To this end, we first show that
for $t\geq1$, the following inequality holds:
\begin{equation}\label{til-v-bdt-bd}
\int_0^1\widetilde{v}^2_t(0,t-x)dx\leq
3\max_{s\in[t-1,t]}\|(\widetilde{v}(\cdot,s),\widetilde{v}_s(\cdot,s))^\top\|_{\mathbb{H}}^2.
\end{equation}
Indeed, define $\rho(t)=2\int_0^1(x-1)\widetilde{v}_t(x,t)
\widetilde{v}_x(x,t)\dd x$. 
Differentiating $\rho(t)$ along the solution of
``$\widetilde{v}$-part'' of \dref{wave-o-hatv} yields
$$\begin{array}{l}
   \dot{\rho}(t)
   \geq\widetilde{v}_t^2(0,t)-\int_0^1[\widetilde{v}
   _x^2(x,t)+\widetilde{v}_t^2(x,t)]dx, \end{array}$$
which, jointly with $|\rho(t)|\leq\|(\widetilde{v}
(\cdot,t),\widetilde{v}_t(\cdot,t))^\top\|_{\mathbb{H}}^2$ and $\int_0^1\widetilde{v}^2_t(0,t-x)\dd x=
\int_{t-1}^t\widetilde{v}_s^2(0,s)\dd s$, implies \dref{til-v-bdt-bd}. Noting that
\begin{equation}\label{W-sol}
   W(x,t)=\left\{\begin{array}{ll}\disp -\frac{c_0}{1-c_0}
   \widetilde{v}(0,t-x),&t\geq x,\crr\disp W_0(x-t), & x>t
   \end{array}\right.
\end{equation}
solves ``$W$-part'' of \dref{wave-o-hatv}. It follows from the
Sobolev embedding theorem and \dref{W-sol} that, for some $C>0$,
\begin{equation}\label{W-sol-est}
\|W(x,t)\|_{H^1(0,1)}\leq C\max_{s\in[t-1,t]}\|(\widetilde{v}(\cdot,s),\widetilde{v}_s(\cdot,s))^\top\|_{\mathbb{H}}^2,
\end{equation}
which gives  $\sup_{t\geq0}\|W(x,t)\|_{H^1(0,1)}<+\infty$.
Since
$\widehat{v}(x,t)=\widetilde{v}(x,t)-W(x,t)$ and $W_t(x,t)=-W_x(x,t)$,
we have $\|(\widehat{v}(\cdot,t),$ $\widehat{v}_t(\cdot,t))^\top\|
_{\mathbb{H}}\leq\|(\widetilde{v}(\cdot,t),\widetilde{v}_t(\cdot,t))
^\top\|_{\mathbb{H}}+\|(W(\cdot,t),$ $W_x(\cdot,t))^\top\|_{\mathbb{H}}$.
Hence, $\sup_{t
\geq0}\|(\widehat{v},\widehat{v}_t,W)(\cdot,t)
^\top\|_{\mathbb{H}\times H^1(0,1)}$ $<+\infty$.

Next, suppose that $\lim_{t\to\infty}|f(w,w_t)|=0$ and $d\in L^2
(0,\infty)$. It follows from Lemma \ref{Lem-ABu} that
``$\widetilde{v}$-part'' of \dref{wave-o-hatv} satisfies $\lim_{t\to
\infty}\|(\widetilde{v}(\cdot,t),\widetilde{v}_t(\cdot,t))^\top\|
_{\mathbb{H}}=0$, 
which, together with \dref{W-sol-est}, leads to $\lim_{t\to\infty}
\|W(x,t)\|_{H^1(0,1)}=0$ and thus 
$\lim_{t\to\infty}\|(\widehat{v}
(\cdot,t),\widehat{v}_t(\cdot,t)\|_{\mathbb{H}}=0$.

Finally, suppose that $f\equiv0$ and $d\equiv0$. Since $\Ascr$ generates
an exponential stable semigroup on $\mathbb{H}$, there exist two
constants $M_2,\mu_2>0$ such that $\|(\widetilde{v}(\cdot,t),
\widetilde{v}_t(\cdot,t))^\top\|_{\mathbb{H}}\leq M_2e^{-\mu_2t}$.
 It follows from \dref{W-sol-est} that $\|W(\cdot,t))\|_{H^1(0,1)}
\leq3CM_2e^{\mu_2}e^{-\mu_2t}$ for all $t\geq
0$, which implies  $\|(\widehat{v}(\cdot,t),\widehat{v}
_t(\cdot,t)\|_{\mathbb{H}}$ $\leq Me^{-\mu t}$ with some $M,\mu>0$.
 \hfill $\Box$

The system \dref{wave-o-v_second} is used to estimate the total
disturbance. Actually, Let $\widetilde{z}(x,t)=z(x,t)-\widehat{v}
(x,t)-W(x,t)$. Then we can see that $\widetilde{z}(x,t)$ is governed
by
\begin{equation} \label{wave-o-tildz}
   \left\{\begin{array}{l} \widetilde{z}_{tt}(x,t) \m=\m
   \widetilde{z}_{xx}(x,t),\crr\disp
   \widetilde{z}_x(0,t) \m=\m \frac{c_1}{1-c_0}\widetilde{z}(0,t)+
   \frac{c_0}{1-c_0}\widetilde{z}_t(0,t),\crr\disp
   \widetilde{z}(1,t)=0. \end{array}\right.
\end{equation}
We consider system \dref{wave-o-tildz} in $\mathbb{H}_0=H^1_R(0,1)
\times L^2(0,1)$, where $H^1_R(0,1)=\{\phi\in H^1(0,1):\phi(1)=0\}$.
Noting that \dref{wave-o-tildz} is exponentially stable on
$\mathbb{H}_0$, the following lemma is easily obtained.

\begin{lemma}\label{dis-est}
Let $c_0\in(0,1)$ and $c_1>0$. For any initial value $(\widetilde{z}_0,\widetilde{z}_1)\in \mathbb{H}_0$, system \dref{wave-o-tildz}
admits a unique solution  $(\widetilde{z},\widetilde{z}_t)\in C(0,\infty;\mathbb{H}_0)$ which satisfies
$\widetilde{z}_x(1,t)\in L^2(0,\infty)$.
\end{lemma}

Since $\widetilde{z}_x(1,t)=z_x(1,t)+F(t)$, where $F$ is as in
\dref{Helen}, by Lemmas \ref{Lem-ABu} and \ref{dis-est} we can regard
$-z_x(1,t)$ as an estimate of the disturbance $F(t)$, that is,
$-z_x(1,t)\approx F(t)$.

\section{Controller and observer design} \label{Sec-stateobserver-U-N}

In this section,  based on our disturbance estimator, we design a
state observer for the system \dref{wave-o} as follows:
\begin{equation}\label{wave-o-obser}
\left\{\begin{array}{l}
\widehat{w}_{tt}(x,t)=\widehat{w}_{xx}(x,t),\;0<x<1,\;t>0, \crr\disp
\widehat{w}_x(0,t)=-qw(0,t)+c_1[\widehat{w}(0,t)-w(0,t)],  \crr\disp
\widehat{w}_x(1,t)=u(t)-z_x(1,t)-Y_x(1,t),\;t\geq0, \crr\disp
\widehat{w}(x,0)=\widehat{w}_0(x),\;\widehat{w}_t(x,0)=\widehat{w}_1
   (x),\; 0\leq x\leq1,\crr\disp
Y_t(x,t)=-Y_x(x,t),\;0<x<1,\;t>0, \crr\disp
Y(0,t)=-c_0[\widehat{w}(0,t)-w(0,t)],\;t\geq0, \crr\disp
Y(x,0)=Y_0(x),\;0\leq x\leq 1,
\end{array}\right.
\end{equation}
where $c_0$ and $c_1$ are the same design parameters as in in
\dref{wave-o-v_first} and \dref{wave-o-v_second}. The signal
$-z_x(1,t)$, generated by the total disturbance estimator, is used to
compensate $F(t)$. The observer \dref{wave-o-obser} is a ``natural
observer'' after canceling the disturbance in a sense that it employs
a copy of the plant plus output injection (in this case, only at the
boundary). Note that the observer \dref{wave-o-obser} is different
from the observer in \cite{KrsGuoetc2008auto}, where the signal
$w_t(1,t)$ (that is considered unavailable in this paper) is used.

To show the asymptotical convergence of the observer above,
we introduce the observer error variable $\varepsilon(x,t)=\widehat{w}
(x,t)-w(x,t)$. Then $(\varepsilon(x,t),Y(x,t))$ satisfies
\begin{equation}\label{wave-o-obser-err-P}
\left\{\begin{array}{l}
\varepsilon_{tt}(x,t)=\varepsilon_{xx}(x,t),  \crr\disp
\varepsilon_x(0,t)\hspace{-0.1cm}=\hspace{-0.1cm}c_1\varepsilon(0,t),\;\;
\varepsilon_x(1,t)\hspace{-0.1cm}=\hspace{-0.1cm}-\widetilde{z}_x(1,t)-Y_x(1,t), \crr\disp
Y_t(x,t)=-Y_x(x,t), \;\; 
Y(0,t)=-c_0\varepsilon(0,t).
\end{array}\right.
\end{equation}
\begin{lemma}\label{lem-varepto0}
Let $c_0\in(0,1)$, $c_1>0$ and let the signal $\widetilde{z}_x(1,t)$ be generated by system \dref{wave-o-tildz}.
Then,  for any initial value
$(\varepsilon(\cdot,0),\varepsilon_t(\cdot,0),Y(\cdot,0))^\top\in \mathbb{H}\times H^1(0,1)$
with the compatibility condition $Y(0,0)=-c_0\varepsilon(0,0)$,
then \dref{wave-o-obser-err-P}  admits  a unique solution  $(\varepsilon,\varepsilon_t,Y)^\top\in
C(0,\infty;\mathbb{H}\times H^1(0,1))$ satisfying
$
\lim_{t\to\infty}\|(\varepsilon(\cdot,t),\varepsilon_t(\cdot,t),Y(\cdot,t))\|_{\mathbb{H}\times H^1(0,1)}=0.
$
\end{lemma}
{\bf Proof.}
We introduce a new variable $\widetilde\varepsilon(x,t)=\varepsilon(x,t)+Y(x,t)$.
 Then,  $(\widetilde\varepsilon(x,t),Y(x,t))$ is governed by
\begin{equation}\label{wave-o-obser-err}
\left\{\begin{array}{l}
\widetilde\varepsilon_{tt}(x,t)=\widetilde\varepsilon_{xx}(x,t),  \crr
\widetilde\varepsilon_x(0,t)=\frac{c_1}{1-c_0}\widetilde\varepsilon(0,t)+\frac{c_0}{1-c_0}\widetilde\varepsilon_t(0,t), \crr\disp
\widetilde\varepsilon_x(1,t)=-\widetilde{z}_x(1,t),  \crr\disp
Y_t(x,t)=-Y_x(x,t),  \;\;
Y(0,t)=-\frac{c_0}{1-c_0}\widetilde\varepsilon(0,t).
\end{array}\right.
\end{equation}
The ``$\widetilde\varepsilon$-part'' of \dref{wave-o-obser-err} can be
$\frac{d}{dt}(\widetilde\varepsilon(\cdot,t),\widetilde\varepsilon_t(\cdot,t))
=\mathcal{A}(\widetilde\varepsilon(\cdot,t),$ $\widetilde\varepsilon_t(\cdot,t))+\mathcal{B}\widetilde{z}_x(1,t)$,
where  $\mathcal{A}$ and $\mathcal{B}$ are defined in the proof of Lemma \ref{Lem-hatz-bd}.
Since ${\mathcal A}$ generates an exponential stable
$C_0$-semigroup $e^{{\mathcal A}t}$ and $\mathcal{B}$ is admissible to
$e^{{\mathcal A}t}$, it follows from Lemma \ref{Lem-ABu} and $\widetilde{z}_x(1,t)\in L^2(0,\infty)$ due to Lemma  \ref{dis-est} that ``$\widetilde\varepsilon$-part'' of \dref{wave-o-obser-err}
has a unique solution that is asymptotically stable.
Next,  $\lim_{t\to\infty}\|Y(\cdot,t))\|_{H^1(0,1)}=0$ can be easily obtained in the same way like Lemma \ref{Lem-hatz-bd}
by noting  the fact that
\begin{equation}\label{Y-wp-sol}
Y(x,t)=
\left\{\begin{array}{ll}
\frac{-c_0}{1-c_0}\widetilde\varepsilon(0,t-x),&t\geq x,\crr
Y_0(x-t),&x>t.
\end{array}\right.
\end{equation}
solves the ``$Y$-part'' of \dref{wave-o-obser-err}.
Since $\varepsilon(x,t)=\widetilde\varepsilon(x,t)-Y(x,t)$
and $Y_t(x,t)=-Y_x(x,t)$,  we have
$
\|(\varepsilon(\cdot,t),\varepsilon_t(\cdot,t))\|_{\mathbb{H}}$ $
\leq \|(\widetilde\varepsilon(\cdot,t),\widetilde\varepsilon_t(\cdot,t))\|_{\mathbb{H}}
\hspace{-0.1cm}+\hspace{-0.1cm}\|(Y(\cdot,t),Y_x(\cdot,t))\|_{\mathbb{H}}\to0,
$ as $t\hspace{-0.1cm}\to\hspace{-0.1cm}\infty$.
$\hspace{0cm}\Box$\\
By Lemma \ref{lem-varepto0}, \dref{wave-o-obser} is a state observer of \dref{wave-o}.
Now, by  the observer-based feedback control law of \cite{KrsGuoetc2008auto},
we propose the following observer-based feedback controller :
\begin{equation}\label{con-out}
\begin{array}{l}\disp
u(t)\hspace{-0.05cm}=\hspace{-0.05cm}z_x(1,t)\hspace{-0.05cm}+\hspace{-0.05cm}Y_x(1,t)\hspace{-0.05cm}-\hspace{-0.05cm}
c_3\widehat{w}_t(1,t)\hspace{-0.05cm}-\hspace{-0.05cm}(c_2\hspace{-0.05cm}+\hspace{-0.05cm}q)\widehat{w}(1,t)\crr\disp\hspace{0.8cm}
-(c_2+q)\int_0^1e^{q(1-\xi)}[c_3\widehat{w}_t(\xi,t)+q\widehat{w}(\xi,t)]d\xi,
\end{array}
\end{equation}
where $c_2,c_3$ are positive design parameters.
The ``$\widehat{w}$-part'' of the closed-loop of observer  \dref{wave-o-obser} corresponding to controller \dref{con-out} becomes
\begin{equation}\label{wave-o-obser-closed}
\left\{\begin{array}{l}
\widehat{w}_{tt}(x,t)=\widehat{w}_{xx}(x,t),  \crr\disp
\widehat{w}_x(0,t)=-qw(0,t)+c_1[\widehat{w}(0,t)-w(0,t)],  \crr\disp
\widehat{w}_x(1,t)=-c_3\widehat{w}_t(1,t)-(c_2+q)\widehat{w}(1,t)\cr\disp\hspace{0.3cm}
-(c_2+q)\int_0^1e^{q(1-\xi)}[c_3\widehat{w}_t(\xi,t)+q\widehat{w}(\xi,t)]d\xi.
\end{array}\right.
\end{equation}
Consider the transformation (\cite{KrsGuoetc2008auto})
\begin{equation}\label{hatw-tw}
\widetilde{w}(x,t)
=\widehat{w}(x,t)+(c_2+q)\int_0^xe^{q(x-\xi)}\widehat{w}(\xi,t)d\xi,
\end{equation}
and its  inverse transformation is given by
\begin{equation}
\widehat{w}(x,t)
=\widetilde{w}(x,t)-(c_2+q)\int_0^xe^{-c_2(x-\xi)}\widetilde{w}(\xi,t)d\xi.
\end{equation}
It can be shown that  \dref{hatw-tw} converts system
\dref{wave-o-obser-closed} into
\begin{equation} \label{tilde-w-sys}
   \left\{\begin{array}{l} \widetilde{w}_{tt}(x,t) =
   \widetilde{w}_{xx}(x,t) - (c_1+q)(c_2+q)e^{qx}\e(0,t),\crr
   \widetilde{w}_x(0,t) \m=\m c_2\widetilde{w}(0,t)+(c_1+q)\e(0,t),\crr
   \widetilde{w}_x(1,t) \m=\m -c_3\widetilde{w}_t(1,t).
   \end{array}\right.
\end{equation}
From Lemma \ref{lem-varepto0}, we can see
$\lim_{t\to\infty}|\varepsilon(0,t)|=0$.  We can show that system
\dref{tilde-w-sys} is asymptotically stable by making use of Lemma
\ref{Lem-ABu}. Actually, we can show that system \dref{tilde-w-sys}
is exponentially stable. For this, we consider the coupled
system consisting of \dref{wave-o-tildz}, \dref{wave-o-obser-err-P}
and \dref{tilde-w-sys} together, in the space $\mathcal{X}=\mathbb{H}
\times H^1(0,1)\times H^1_R(0,1)\times L^2(0,1)\times \mathbb{H}$.

\begin{theorem}\label{calADAto0}
Suppose that $c_0\in(0,1)$ and $c_i>0$, $i=1,2,3$. For any initial
value $(\widetilde\e_0,\widetilde\e_1,Y_0,\widetilde{z}_0,\widetilde
{z}_t,\widetilde{w}_0,\widetilde{w}_1)\in\mathcal{X}$, with the
compatibility condition $Y_0(0)=-c_0\widetilde\e_0(0)$, there exists
a unique solution $(\e,\e_t,Y,\widetilde{z},\widetilde{z}_t,
\widetilde{w},\widetilde{w}_t)\in C(0,\infty;\mathcal{X})$ to
\dref{wave-o-tildz}, \dref{wave-o-obser-err-P} and \dref{tilde-w-sys}
such that
$ \|(\e,\e_t,Y,\widetilde{z},\widetilde{z}_t,\widetilde{w},
   \widetilde{w}_t)(\cdot,t)\|_{\mathcal{X}}$ $ \leq Me^{-\mu t}$
with some $M,\mu>0$.
\end{theorem}

{\bf Proof.} Let  $\widetilde\e(x,t)=\varepsilon(x,t)+Y(x,t)$ and
$\eta(x,t)=\widetilde\e(x,t)+\widetilde{z}(x,t)$. Then, $(\eta(x,t),
Y(x,t))$ satisfies the following PDEs:
\begin{equation} \label{wave-o-obser-err-and-tilz-tilw-equiv}
   \left\{\begin{array}{l} \eta_{tt}(x,t) \m=\m \eta_{xx}(x,t),\cr
   \eta_x(0,t) \m=\m \frac{c_1}{1-c_0}\eta(0,t)+\frac{c_0}{1-c_0}
   \eta_t(0,t),\crr \eta_x(1,t) \m=\m 0,
   Y_t(x,t) \m=\m -Y_x(x,t),
   \crr Y(0,t) \m=\m -\frac{c_0}{1-c_0}[\eta(0,t)-\widetilde{z}(0,t)].
\end{array}\right.
\end{equation}
It is well known that the ``$\eta$-part'' of
\dref{wave-o-obser-err-and-tilz-tilw-equiv} and also \dref{wave-o-tildz}
are exponentially stable, which implies that $(\widetilde\e,\widetilde
\e_t)$ is also exponentially stable on $\mathbb{H}$. 
Similar to the proof of the last assertion of Lemma \ref{Lem-hatz-bd},
we obtain the exponential stability of ``$Y$-part'' of
\dref{wave-o-obser-err-and-tilz-tilw-equiv}. Since $\e=\widetilde\e-Y$
and $\e_t=\widetilde\e_t+Y_x$, we have that $(\e,\e_t)$ is exponentially
stable on $\mathbb{H}$. By the Sobolev embedding theorem, $|\e(0,t)|\leq
\|\e\|_{H^1(0,1)}$, which shows that $\e(0,t)$ decays exponentially.
Rewrite \dref{tilde-w-sys} as
\begin{equation} \label{varep-abs-tep}
   \frac{\dd}{\dd t} \begin{pmatrix}\widetilde{w}(\cdot,t)\cr
   \widetilde{w}_t(\cdot,t)\end{pmatrix} = A_{0}\begin{pmatrix}
   \widetilde{w}(\cdot,t)\cr \widetilde{w}_t(\cdot,t)\end{pmatrix}
   +{B}_1\e(0,t)+{B}_2\e(0,t),
\end{equation}
where the operator $A_{0}:D(A_{0})(\subset\mathbb{H})\to \mathbb{H}$ is
given by
\begin{equation} \label{Aw-def}
   \left\{\begin{array}{l}\disp A_{0}(\phi,\psi)^\top \m=\m (\psi,
   \phi'')^\top\ \ \forall (\phi,\psi)^\top\in D({A}_{0}),\crr\disp
   D(A_0) \m=\m \big\{(\phi,\psi)^\top\in H^2(0,1)\times H^1(0,1):\crr
   \quad\quad\phi'(0) \m=\m c_2\phi(0),\ \ \phi'(1) \m=\m -c_3\psi(1)\big\},
   \end{array}\right.
\end{equation}
and $B_1=(c_1+q)(0,-\delta(x))$, $B_2=-(c_1+q)(c_2+q)(0,-e^{qx})$.
Since $A_{0}$ generates an exponentially stable $C_0$-semigroup
$e^{A_{0}t}$ on $\mathbb{H}$ and $B_1$, $B_2$ are admissible for
$e^{A_{0}t}$. By Lemma \ref{Lem-ABu}, we know that $(\widetilde{w},
\widetilde{w}_t)$ is exponentially stable on $\mathbb{H}$. \hfill
$\hspace{3cm}\Box$

\section{Well-posedness and exponential stability of the closed-loop
         system} \label{Sec-closed-U-N} 

We go back to the closed-loop system \dref{wave-o} under the feedback
\dref{con-out}:
\begin{equation}\label{wave-o-closed}
\left\{\begin{array}{l}
w_{tt}(x,t)=w_{xx}(x,t),  \crr\disp
w_x(0,t)=-qw(0,t), \crr\disp
w_x(1,t)=z_x(1,t)+Y_x(1,t)-(c_3+q)\widehat{w}(1,t)\crr\disp\hspace{0.3cm}
      -c_3\widehat{w}_t(1,t)+f(w(\cdot,t),w_t(\cdot,t))+d(t) \crr\disp\hspace{0.3cm}
     -(c_3+q)\int_0^1e^{q(1-\xi)}[c_3\widehat{w}_t(\xi,t)+q\widehat{w}(\xi,t)]d\xi, \crr\disp
v_{tt}(x,t)=v_{xx}(x,t),   \crr\disp
v_x(0,t)=-qw(0,t)+c_1[v(0,t)-w(0,t)], \crr\disp
v_x(1,t)=z_x(1,t)+Y_x(1,t)-W_x(1,t)\crr\disp
\hspace{0.2cm}-c_3\widehat{w}_t(1,t)-(c_2+q)\widehat{w}(1,t)
\crr\disp\hspace{0.2cm}-(c_2+q)\int_0^1e^{q(1-\xi)}[c_3\widehat{w}_t(\xi,t)+q\widehat{w}(\xi,t)]d\xi, \crr\disp
z_{tt}(x,t)=z_{xx}(x,t),  \crr\disp
z_x(0,t)=\frac{c_1}{1-c_0}z(0,t)+\frac{c_0}{1-c_0}z_t(0,t),  \crr\disp
z(1,t)=v(1,t)+W(1,t)-w(1,t),  \crr\disp
W_t(x,t)=-W_x(x,t), \crr\disp
W(0,t)=-c_0[v(0,t)-w(0,t)], \crr\disp
\widehat{w}_{tt}(x,t)=\widehat{w}_{xx}(x,t),   \crr\disp
\widehat{w}_x(0,t)=-qw(0,t)+c_1[\widehat{w}(0,t)-w(0,t)],  \crr\disp
\widehat{w}_x(1,t)=-c_3\widehat{w}_t(1,t)-(c_2+q)\widehat{w}(1,t)\crr\disp\hspace{0.3cm}
-(c_2+q)\int_0^1e^{q(1-\xi)}[c_3\widehat{w}_t(\xi,t)+q\widehat{w}(\xi,t)]d\xi, \crr\disp
Y_t(x,t)=-Y_x(x,t),   \crr\disp
Y(0,t)=-c_0[\widehat{w}(0,t)-w(0,t)].  
\end{array}\right.
\end{equation}
We consider system \dref{wave-o-closed} in the state space $\mathcal{H}=\mathbb{H}^3\times H^1(0,1)\times \mathbb{H}\times H^1(0,1)$.
\begin{theorem}
Suppose that $c_0\in(0,1)$ and $c_i>0$, $i=1,2,3$, $f:\mathbb{H}\to\mathbb{R}$ is continuous, and $d\in L^\infty(0,\infty)$ or $d\in L^2(0,\infty)$.
For any initial value $(w_0,w_1,v_0,v_1,z_0,z_1,W_0,\widehat{w}_0,$ $\widehat{w}_1,Y_0)\in \mathcal{H}$ with the compatibility
conditions
$
z_0(1)-v_0(1)-W_0(1)+w_0(1)=0,\;W_0(0)+c_0[v_0(0)-w_0(0)]=0,\; Y_0(0)+c_0[\widehat{w}_0(0)-w_0(0)]=0,
$
then \dref{wave-o-closed} admits a unique solution  $(w,w_t,v,v_t,W,z,z_t,\widehat{w}_0,\widehat{w}_t,Y)\in C(0,\infty;\mathcal{H})$ satisfying
$
\|(w(\cdot,t),w_t(\cdot,t), \widehat{w}(\cdot,t),\widehat{w}_t(\cdot,t),Y(\cdot,t))
\|_{\mathbb{H}^2\times H^1(0,1)}$ $\leq Me^{-\mu t},\forall t\geq0,
$
with some $M,\mu>0$, and
$
\sup_{t\geq0}\|(v(\cdot,t),$ $v_t(\cdot,t),z(\cdot,t),z_t(\cdot,t),W(\cdot,t))\|_{\mathbb{H}^2\times H^1(0,1)}<+\infty.
$
Moreover,  if $f(0,0)=0$ and $d\in L^2(0,\infty)$, then $\lim_{t\to\infty}\|(v(\cdot,t),$ $v_t(\cdot,t),z(\cdot,t),
z_t(\cdot,t),W(\cdot,t))\|_{\mathbb{H}^2\times H^1(0,1)}=0$. If $f\equiv0$ and $d\equiv0$, then
$\|\hspace{-0.03cm}(v(\cdot,t),\hspace{-0.03cm}v_t(\cdot,t),\hspace{-0.03cm}z(\cdot,t),
\hspace{-0.03cm}z_t(\cdot,t),\hspace{-0.03cm}W(\cdot,t))
\hspace{-0.02cm}\|_{\mathbb{H}^2\hspace{-0.03cm}\times\hspace{-0.03cm} H^1(0,1)}$ $\leq M'e^{-\mu' t}$ for all $t\geq0$
with some $M',\mu'>0$.
\end{theorem}
Using Lemma \ref{Lem-hatz-bd} and Theorem \ref{calADAto0}, the above theorem can be established. We skip the details for lack of space.

\section{The anti-stable wave equation} \label{Sec-US-N} 

In this section we consider the output feedback exponential
stabilization for the system governed by the following equations,
where $t>0$:
\begin{equation}\label{wave-o-US-N}
\left\{\begin{array}{l}
   w_{tt}(x,t)=w_{xx}(x,t),\;0<x<1,\crr
   w_x(0,t)=-qw_t(0,t),\crr
   w_x(1,t)=u(t) + f(w(\cdot,t),w_t(\cdot,t)) + d(t),\crr
   w(x,0) = w_0(x),\ \ w_t(x,0) = w_1(x),\ \ 0\leq x\leq 1,\crr
   y_{m}(t) = \{w(0,t),w(1,t)\},\end{array}\right.
\end{equation}
where $(w,w_t)$ is the state, $u$ is the control input signal, $y_m$
is the output signal and $q>0$ with $q\neq1$. As in \dref{wave-o},
$f:H^1(0,1)\times L^2(0,1)\to\rline$ is an unknown possibly nonlinear
mapping that represents the internal uncertainty, and $d\in L^\infty
(0,\infty)$ represents the unknown external disturbance.

\subsection{Disturbance estimator}
We design, in terms of $y_{m}(t)=\{w(0,t),w(1,t)\}$, a disturbance
estimator for system \dref{wave-o-US-N} as follows: for all $t>0$,
\begin{equation} \label{wave-o-v-US-N}
   \left\{\begin{array}{l} v_{tt}(x,t) = v_{xx}(x,t),\ \ 0<x<1,\crr
   v_x(0,t) = -qv_t(0,t) + c_1[v(0,t)-w(0,t)], \crr
   v_x(1,t) = u(t)-W_x(1,t),\crr
   z_{tt}(x,t) = z_{xx}(x,t),\ \ 0<x<1,\cr
   z_x(0,t) = \frac{c_1}{1-c_0}z(0,t) +
      \frac{c_0-q}{1-c_0}z_t(0,t),\cr
   z(1,t) = v(1,t)+W(1,t)-w(1,t),\crr
   W_t(x,t) = -W_x(x,t),\;0<x<1,\crr
   W(0,t) = -c_0[v(0,t)-w(0,t)], \end{array}\right.
\end{equation}
where $c_0$ and $c_1$ are two design parameters so that
$\frac{c_1}{1-c_0}>0$ and $\frac{c_0-q}{1-c_0}>0$.  The above
disturbance estimator \dref{wave-o-v-US-N} is uniquely determined by
input signal $u$ and two measurement signals $\{w(0,t),w(1,t)\}$. In
the disturbance estimator, the ``$(v,W)$-part'' is used to channel total
disturbance from original system. Indeed, let $\widehat{v}(x,t)=v(x,t)
-w(x,t)$. Then it is easy to check that $\widehat{v}(x,t)$ satisfies
\begin{equation} \label{wave-o-hatv-P-US-N}
\left\{\begin{array}{l}
\widehat{v}_{tt}(x,t)=\widehat{v}_{xx}(x,t),  \crr\disp
\widehat{v}_x(0,t)=-q\widehat{v}_t(0,t)+c_1\widehat{v}(0,t),  \crr\disp
\widehat{v}_x(1,t)=-f(w(\cdot,t),w_t(\cdot,t))-d(t)-W_x(1,t), \crr\disp
W_t(x,t)=-W_x(x,t),\;\;
W(0,t)=-c_0\widehat{v}(0,t).
\end{array}\right.
\end{equation}
Similar to Lemma \ref{Lem-hatz-bd}, we can prove the following:

\begin{lemma}\label{Lem-hatz-bd-US-N}
Suppose that $\frac{c_1}{1-c_0}>0$, $\frac{c_0-q}{1-c_0}>0$; $d\in
L^\infty(0,\infty)$, (or $d\in L^2(0,\infty)$),
$f:\mathbb{H}\to\mathbb{R}$ is continuous and that \dref{wave-o-US-N}
admits a unique solution $(w,\dot{w})^\top\in C(0,\infty;\mathbb{H})$
which is bounded. For any initial value
$(\widehat{v}_0,\widehat{v}_1,W_0)^\top\in \mathbb{H}\times H^1(0,1)$
with the compatibility condition $W_0(0)=-c_0\widehat{v}_0(0)$, system
\dref{wave-o-hatv-P-US-N} admits a unique solution
$(\widehat{v},\widehat{v}_t,W)^\top\in C(0,\infty;\mathbb{H}\times
H^1(0,1))$ such that $\sup_{t\geq0}\|(\widehat{v},
\widehat{v}_t,W)(\cdot,t))^\top\|_{\mathbb{H}\times
H^1(0,1)}<+\infty$. Moreover, if $f(0,0)=0$ and $d\in L^2(0,\infty)$,
then $\lim_{t\to\infty}\|(\widehat{v},\widehat{v}_t
,W)(\cdot,t)^\top\|_{\mathbb{H}\times H^1(0,1)}=0$. If
$f\equiv d\equiv0$, then $\|(\widehat{v},
\widehat{v}_t,W)(\cdot,t)^\top\|_{\mathbb{H}\times
H^1(0,1)}\leq M'e^{-\mu' t}$ for all $t\geq0$ with some $M',\mu'>0$.
\end{lemma}

Let $\widetilde{z}(x,t)\hspace{-0.05cm}=\hspace{-0.05cm}z(x,t)
\hspace{-0.05cm}-\hspace{-0.05cm}\widehat{v}(x,t)
\hspace{-0.05cm}-\hspace{-0.05cm}W(x,t)$. Then $\widetilde{z}(x,t)$
satisfies
\begin{equation} \label{wave-o-tildz-US-N}
   \left\{\begin{array}{l} \widetilde{z}_{tt}(x,t) = \widetilde{z}
   _{xx}(x,t),\crr \widetilde{z}_x(0,t) \hspace{-0.1cm}=\hspace{-0.1cm} \frac{c_1}{1-c_0}
   \widetilde{z}(0,t)\hspace{-0.1cm}+\hspace{-0.1cm}\frac{c_0-q}{1-c_0}\widetilde{z}_t(0,t),\;\;
   \widetilde{z}(1,t)\hspace{-0.1cm}=\hspace{-0.1cm}0, \end{array}\right.
\end{equation}
which is exactly the same as the system \dref{wave-o-tildz} by
replacing $\frac{c_0}{1-c_0}$ with $\frac{c_0-q}{1-c_0}$.  Thus, it
follows from Lemma \ref{dis-est} that $-z_x(1,t)$ can be regarded as
an estimation of the disturbance $f(w(\cdot,t),w_t(\cdot,t))+d(t)$.

\subsection{Controller and observer design}

Using the disturbance estimator \dref{wave-o-v-US-N}, we now design
the following observer for system \dref{wave-o-US-N}:
\begin{equation} \label{wave-o-obser-US-N}
   \left\{\begin{array}{l} \widehat{w}_{tt}(x,t) = \widehat{w}_{xx}
   (x,t),\;0<x<1,\;t>0, \crr \widehat{w}_x(0,t) = -q\widehat{w}_t(0,t)
   +c_1[\widehat{w}(0,t)-w(0,t)],\crr \widehat{w}_x(1,t) = u(t)-z_x
   (1,t)-Y_x(1,t),\;t\geq 0,\crr \widehat{w}(x,0) = \widehat{w}_0(x),
   \ \ \widehat{w}_t(x,0) = \widehat{w}_1(x),\crr
   Y_t(x,t) = -Y_x(x,t),\;0<x<1,\;t>0,\crr
   Y(0,t) = -c_0[\widehat{w}(0,t)-w(0,t)],\;t\geq 0,\crr
   Y(x,0) = Y_0(x),\; 0\leq x\leq 1, \end{array}\right.
\end{equation}
where $c_1$ and $c_2$ are two design parameters that are the same as
in \dref{wave-o-v-US-N}; $-z_x(1,t)$ plays the role of total
disturbance $f(w(\cdot,t),w_t(\cdot,t))+d(t)$. Let $\e(x,t)=\widehat{w}
(x,t)-w(x,t)$ be the observer error. Then $(\e(x,t),Y(x,t))$ is
governed by
\begin{equation} \label{wave-o-obser-err-P-US-N}
\left\{\begin{array}{l} \varepsilon_{tt}(x,t)=\varepsilon_{xx}(x,t),
\crr\disp \varepsilon_x(0,t)=-q\varepsilon_t(0,t)+c_1\varepsilon(0,t),
\crr\disp \varepsilon_x(1,t)=-\widetilde{z}_x(1,t)-Y_x(1,t), \crr\disp
Y_t(x,t)=-Y_x(x,t), \;\;
Y(0,t)=-c_0\varepsilon(0,t).
\end{array}\right.
\end{equation}
Similar to Lemma \ref{lem-varepto0},   introducing a new variable $\widetilde\varepsilon(x,t)=\varepsilon(x,t)+Y(x,t)$,
we can establish the following Lemma:
\begin{lemma}\label{lem-varepto0-US-N}
Suppose that $\frac{c_0-q}{1-c_0}>0$, $\frac{c_1}{1-c_0}>0$ and
the signal $\widetilde{z}_x(1,t)$ is generated  by system \dref{wave-o-tildz}.
Then,  for any initial value $(\varepsilon(\cdot,0),\varepsilon_t(\cdot,0),Y(\cdot,0))^\top\in \mathbb{H}\times H^1(0,1)$
with the compatibility condition $Y(0,0)=-c_0\varepsilon(0,0)$,
then \dref{wave-o-obser-err-P-US-N} admits  a unique solution  $(\varepsilon,\varepsilon_t,Y)^\top\in
C(0,\infty;\mathbb{H}\times H^1(0,1))$ satisfying
$
\lim_{t\to\infty}\|(\varepsilon(\cdot,t),\varepsilon_t(\cdot,t),Y(\cdot,t))\|_{\mathbb{H}\times H^1(0,1)}=0.
$
\end{lemma}
By Lemma \ref{lem-varepto0-US-N},
to find a stabilizing control law for system \dref{wave-o-US-N},
it suffices to find a stabilizing control law for system \dref{wave-o-obser-US-N}. To this end,
we first introduce an auxiliary system:
\begin{equation}\label{Z-aux-sys}
\left\{\begin{array}{l}
Z_t(x,t)=-Z_x(x,t),\;0<x<1,\;t>0, \crr\disp
Z(0,t)=-c_2\widehat{w}(0,t),\;t\geq0, \crr\disp
Z(x,0)=Z_0(x),\;0\leq x\leq1.
\end{array}\right.
\end{equation}
Since the above system depends on the observer \dref{wave-o-obser-US-N},
it is implementable.
Next, we introduce $\widetilde{w}(x,t)=\widehat{w}(x,t)+Z(x,t)$. Then, $(\widetilde{w}(x,t),Z(x,t))$ satisfies
\begin{equation}\label{wave-o-obser-trans-tw-US-N}
\left\{\begin{array}{l}
\widetilde{w}_{tt}(x,t)=\widetilde{w}_{xx}(x,t),  \crr\disp
\widetilde{w}_x(0,t)=\frac{c_2-q}{1-c_2}\widetilde{w}_t(0,t)+c_1\varepsilon(0,t),  \crr\disp
\widetilde{w}_x(1,t)=u(t)-z_x(1,t)-Y_x(1,t)+Z_x(1,t),  \crr\disp
Z_t(x,t)=-Z_x(x,t), \;\; 
Z(0,t)=-\frac{c_2}{1-c_2}\widetilde{w}(0,t). 
\end{array}\right.
\end{equation}
It is seen that there is a ``passive damper'' at the left end $x=0$ if $\frac{c_2-q}{1-c_2}>0$ and $\varepsilon(0,t)=0$.
It is seen that the exponential stability of system \dref{wave-o-obser-trans-tw-US-N}
 is equivalent to the exponential stability of system \dref{wave-o-obser-US-N}.
We propose the following observer-based feedback controller as follows:
\begin{equation}\label{con-out-US-N}
\begin{array}{l}\disp
u(t)=-c_3\widehat{w}(1,t)-c_3Z(1,t)+z_x(1,t)\crr\disp\hspace{1cm}+Y_x(1,t)-Z_x(1,t).
\end{array}
\end{equation}
Under the control law \dref{con-out-US-N}, the ``$\widetilde{w}$-part'' of the closed-loop of  \dref{wave-o-obser-trans-tw-US-N}  becomes
\begin{equation}\label{wave-o-obser-trans-tw-closed-US-N}
\left\{\begin{array}{l}
\widetilde{w}_{tt}(x,t)=\widetilde{w}_{xx}(x,t),  \crr\disp
\widetilde{w}_x(0,t)=\frac{c_2-q}{1-c_2}\widetilde{w}_t(0,t)+c_1\varepsilon(0,t),  \crr\disp
\widetilde{w}_x(1,t)=-c_3\widetilde{w}(1,t),  \crr\disp
Z_t(x,t)=-Z_x(x,t),  \;\;
Z(0,t)=-\frac{c_2}{1-c_2}\widetilde{w}(0,t).
\end{array}\right.
\end{equation}
Using Lemma \ref{lem-varepto0-US-N}, we can obtain the asymptotical
stability of the solution of \dref{wave-o-obser-trans-tw-closed-US-N}.
Actually, we can show that system \dref{wave-o-obser-trans-tw-closed-US-N} is exponentially stable.
For this purpose, we consider the coupled system consisting of \dref{wave-o-tildz-US-N}, \dref{wave-o-obser-err-P-US-N},
 \dref{Z-aux-sys} and \dref{wave-o-obser-trans-tw-closed-US-N}
in the space $\mathscr{X}=\mathbb{H}\times H^1(0,1)\times H^1_R(0,1)\times L^2(0,1)\times \mathbb{H}\times H^1(0,1)$.
By  mimicking the proof of Theorem \ref{calADAto0}, we can obtain:
\begin{theorem}\label{calADAto0-US-N}
Suppose that $\frac{c_1}{1-c_0}>0$, $\frac{c_0-q}{1-c_0}>0$, $\frac{c_2-q}{1-c_2}>0$ and $c_3>0$. For any initial value $(\widetilde\varepsilon_0,\widetilde\varepsilon_1,Y_0,\widetilde{z}_0,\widetilde{z}_t,$ $\widetilde{w}_0,\widetilde{w}_1,Z)\in\mathscr{X}$,
with the compatibility conditions $Y_0(0)=-c_0\widetilde\varepsilon_0(0)$, $Z_0(0)=-\frac{c_2}{1-c_2}\widetilde{w}_0(0)$,
then there exists a unique solution
 $(\varepsilon,\varepsilon_t,Y,\widetilde{z},\widetilde{z}_t,\widetilde{w},\widetilde{w}_t,Z)\in C(0,\infty;\mathscr{X})$
to the coupled system consisting of \dref{wave-o-tildz-US-N}, \dref{wave-o-obser-err-P-US-N} and \dref{wave-o-obser-trans-tw-closed-US-N} such that
$
\|(\varepsilon,\varepsilon_t,Y,
\widetilde{z},\widetilde{z}_t,
\widetilde{w},\widetilde{w}_t,Z)(\cdot,t))\|_{\mathscr{X}} \leq Me^{-\mu t}$
 with some $M,\mu>0$.
\end{theorem}

\subsection{Well-posedness and exponential stability of closed-loop system}

We go back to the closed-loop system \dref{wave-o-US-N} under the
feedback \dref{con-out-US-N}:
\begin{equation}\label{wave-o-closed-US-N}
\left\{\begin{array}{l}
w_{tt}(x,t)=w_{xx}(x,t),  \crr\disp
w_x(0,t)=-qw_t(0,t),  \crr\disp
w_x(1,t)=-c_3\widehat{w}(1,t)-c_3Z(1,t)+z_x(1,t)\crr\disp\hspace{0cm}
         +Y_x(1,t)-Z_x(1,t)+f(w(\cdot,t),w_t(\cdot,t))\hspace{-0.05cm}+\hspace{-0.05cm}d(t),  \crr\disp
v_{tt}(x,t)=v_{xx}(x,t),   \crr\disp
v_x(0,t)=-qv_t(0,t)+c_1[v(0,t)-w(0,t)], \crr\disp
v_x(1,t)=-c_3\widehat{w}(1,t)-c_3Z(1,t)+z_x(1,t)\crr\disp\hspace{1cm}
         +Y_x(1,t)-Z_x(1,t)-W_x(1,t),  \crr\disp
z_{tt}(x,t)=z_{xx}(x,t),\;  \crr\disp
z_x(0,t)=\frac{c_1}{1-c_0}z(0,t)+\frac{c_0-q}{1-c_0}z_t(0,t), \crr\disp
z(1,t)=v(1,t)+W(1,t)-w(1,t),  \crr\disp
\widehat{w}_{tt}(x,t)=\widehat{w}_{xx}(x,t),  \crr\disp
\widehat{w}_x(0,t)=-q\widehat{w}_t(0,t)+c_1[\widehat{w}(0,t)-w(0,t)], \crr\disp
\widehat{w}_x(1,t)=-c_3\widehat{w}(1,t)-c_3Z(1,t)-Z_x(1,t), \crr\disp
\hspace{-0.1cm}W_t(x,t)\hspace{-0.1cm}=\hspace{-0.1cm}-\hspace{-0.1cm}W_x(x,t), 
W(0,t)\hspace{-0.1cm}=\hspace{-0.1cm}-c_0[v(0,t)\hspace{-0.1cm}-\hspace{-0.1cm}w(0,t)], \crr\disp
Y_t(x,t)\hspace{-0.1cm}=\hspace{-0.1cm}Y_x(x,t), \;\;
Y(0,t)\hspace{-0.1cm}=\hspace{-0.1cm}-c_0[\widehat{w}(0,t)\hspace{-0.1cm}-\hspace{-0.1cm}w(0,t)],  \crr\disp
Z_t(x,t)=-Z_x(x,t),\;\;
Z(0,t)=-c_2\widehat{w}(0,t).
\end{array}\right.
\end{equation}
We consider system \dref{wave-o-closed-US-N} in the state space
$\mathscr{H}=\mathbb{H}^3\times H^1(0,1)\times \mathbb{H}\times [H^1(0,1)]^2$.
\begin{theorem}\label{Thm-closed-USN}
Suppose that $\frac{c_1}{1-c_0}>0$, $\frac{c_0-q}{1-c_0}>0$, $\frac{c_2-q}{1-c_2}>0$ and $c_3>0$.
Suppose that $f:\mathbb{H}\to\mathbb{R}$ is continuous, and $d\in L^\infty(0,\infty)$ or $d\in L^2(0,\infty)$.
For any initial value $(w_0,w_1,v_0,v_1,z_0,z_1,W_0,\widehat{w}_0,\widehat{w}_1,Y_0)\in \mathscr{H}$ with  the compatibility
condition
$z_0(1)-v_0(1)-W_0(1)+w_0(1)=0,$\\ $Z_0(0)+c_2\widehat{w}_0(0)=0,
W_0(0)+c_0[v_0(0)-w_0(0)]=0,
Y_0(0)+c_0[\widehat{w}_0(0)-w_0(0)]=0,
$
then there exists a unique solution to \dref{wave-o-closed-US-N} such that
 $(w,w_t,v,v_t,W,z,z_t,\widehat{w},\widehat{w}_t,Y)\in C(0,\infty;\mathscr{H})$ satisfying
$
\|(w,\hspace{-0.03cm}w_t,\hspace{-0.03cm}\widehat{w},\hspace{-0.03cm}
\widehat{w}_t,\hspace{-0.03cm}Y,\hspace{-0.03cm}Z)(\cdot,t)$ $
\|_{\mathbb{H}^2\hspace{-0.05cm}\times\hspace{-0.05cm}[H^1(0,1)]^2}$ $\leq Me^{-\mu t},
\forall t\geq0$,
with some $M,\mu>0$, and
$
\sup_{t\geq0}\|(v,v_t,z,$ $z_t,W)(\cdot,t)\|_{\mathbb{H}^2\times H^1(0,1)}$ $<+\infty.
$
Moreover,  if $f(0,0)=0$ and $d\in L^2(0,\infty)$, then $\lim_{t\to\infty}\|(v,
v_t,z,z_t,W)(\cdot,t)\|_{\mathbb{H}^2\times H^1(0,1)}$ $=0$.
 If $f\equiv d\equiv0$, then
$\|(v,
v_t,z,z_t,W)(\cdot,t)\|_{\mathbb{H}^2\times H^1(0,1)}\leq M'e^{-\mu' t}$ for all $t\geq0$
with some $M',\mu'>0$.
\end{theorem}
The above theorem can be established by making use of Lemma \ref{Lem-hatz-bd-US-N} and Theorem \ref{calADAto0-US-N}. We omit
the details owing to space constraints.

\section{Conclusion}

In this paper, we have studied the exponential stabilization problem for 
one dimensional unstable/anti-stable wave equation with Neumann boundary control
 subject to the general disturbance  via the two displacement signals.
The main contributions of this paper are (i) the closed-loop system is exponentially 
stable; (ii) we  use only two output signals which are almost the minimal measurement signals.
The idea used here is potentially promising for treating a moving boundary system, like wave 
PDE dynamics with a moving controlled boundary.




\bibliography{ifacconf}             

\end{document}